\documentclass[reqno,11pt]{amsart}%
\usepackage{amsfonts}
\usepackage{amsmath}
\usepackage{amssymb}
\usepackage{graphicx}%
\providecommand{\U}[1]{\protect\rule{.1in}{.1in}}
\theoremstyle{plain}

\newtheorem{example}{Example}

\numberwithin{equation}{section}
\makeatletter
\@namedef{subjclassname@2020}{\textup{2020} Mathematics Subject Classification}
\makeatother
\begin{document}
\title[Anti-Invariant Holomorphic Statistical Submersions ]{Anti-Invariant Holomorphic Statistical Submersions }
\author{Sema Kazan}
\curraddr{Department of Mathematics, Faculty of Arts and Sciences, Inonu University,
Malatya, Turkey}
\email{sema.bulut@inonu.edu.tr}
\author{Kazuhiko Takano}
\curraddr{Department of Mathematics, School of General Education, Shinshu University,
Nagano 390-8621, Japan}
\email{ktakano@shinshu-u.ac.jp}
\thanks{}
\date{2022}

\begin{abstract}
Our purpose in this article is to study anti-invariant statistical submersions
from holomorphic statistical manifolds. Firstly we introduce holomorphic
statistical submersions satisfying the certain condition, after we give
anti-invariant statistical submersions satisfying the certain condition. And
we supported our results with examples.

\end{abstract}
\maketitle

\section{\textbf{Introduction}}

\bigskip

In 1945, the theory of statistical manifolds has started with a paper of C.R.
Rao \cite{Rao}.

It is known that the theory of statistical manifolds is called as information
geometry. The information geometry, which is typically deals with the study of
various geometric structures on a statistical manifold, has begun as a study
of the geometric structures possessed by a statistical model of probability
distributions. Nowadays, the information geometry has an important application
area, such as, information theory, stochastic processes, dynamical systems and
times series, statistical physics, quantum systems and the mathematical theory
of neural networks \cite{Ay}. Also, some applications of statistical manifolds
in information geometry have been handled in many studies. In \cite{Ali}, the
authors have presented an analytical computation of the asymptotic temporal
behavior of the information geometric complexity of finite dimensional
Gaussian statistical manifolds in the presence of microcorrelations
(correlations between microvariables) and in \cite{Gomez}, the author has
presented an extension of the ergodic, mixing and Bernoulli levels of the
ergodic hierarchy for statistical models on curved manifolds, making use of
elements of the information geometry.

The notion of dual connection (or conjugate connection) in affine geometry,
has been first introduced into statistics by S. Amari \cite{Amari} in 1985. A
statistical model equipped with a Riemannian metric together with a pair of
dual affine connections is called a \textit{statistical manifold}. For more
information about statistical manifolds and information geometry, we refer to
\cite{Blaga}, \cite{Furuhata2}, \cite{Matsuzoe}, \cite{Noguchi} \cite{Calin},
\cite{Kurose} and etc.

Considering these notions, the differential geometry of statistical manifolds
are being studying by geometers by adding different geometric structures to
these manifolds. For instance, in \cite{Vilcu} quaternionic K\"{a}hler-like
statistical manifold have been studied and in \cite{Furuhata}, the authors
have introduced the notion of Sasakian statistical structure and obtained the
condition for a real hypersurface in a holomorphic statistical manifold to
admit such a structure. In \cite{Ahmet}, the author has studied
conformally-projectively flat trans-Sasakian statistical manifolds. Also, the
authors have examined Sasakian statistical manifolds with semi-symmetric
metric connection in \cite{kazanlar}.

Nowadays, some authors has studied statistical submersions. The notion of
statistical submersion between statistical manifolds has introduced in 2001 by
N. Abe and K. Hasegawa \cite{Abe}, the authors generalizing some basic results
of B. O'Neill (\cite{Oneill},\cite{Oneill 3}) concerning Riemannian
submersions and geodesics. Later, K.Takano has introduced statistical
manifolds with almost complex structures and its submersions \cite{takano} in
2004. Also, in \cite{Takano 4}, Takano has given examples of the statistical
submersion and in \cite{Takano 2} has studied statistical submersions of
statistical manifolds with almost contact structures. Quaternionic
K\"{a}hler-like statistical submersions has been given in \cite{Vilcu}. In
\cite{Vilcu 2}, G.E. Vilcu has studied para-K\"{a}hler-like statistical
submersions. For other works see \cite{Aliya}, \cite{Hulya}. Later such
submersions have been considered between manifolds with differentiable
structures. B. Watson defined almost Hermitian submersions between almost
Hermitian manifolds and he showed that the base manifold and each fibre have
the same kind of structure as the total space, in most cases \cite{Watson}.
And, many authors have studied on submersions, see \cite{Sahin}, \cite{Sahin
2}, \cite{Sahin 3}, \cite{Johnson}, \cite{Lee} \cite{Akif}, \cite{Akif 3},
\cite{Akif 2}, \cite{Akif 4}.

In Sect.2, we introduce a brief introduction about statistical manifolds and
we give the definition and example of the holomorphic statistical manifolds.
In Sect.3, we investigate holomorphic statistical submersions satisfying the
certain condition. We give an example of holomorphic statistical submersion.
In Sect.4, we define the anti-invariant statistical submersion from
holomorphic statistical manifolds and we study anti-invariant statistical
submersions satisfying the certain conditions. We give an example and some results.

\section{\textbf{Holomorphic Statistical Manifolds}}

An $m$-dimensional semi-Riemannian manifold is a smooth manifold
$\mathcal{M}^{m}$ furnished with a metric $g$, where $g$ is a symmetric
nondegenerate tensor field on $\mathcal{M}$ of constant index. The common
value $\nu$ of index $g$ on $\mathcal{M}$ is called the index of
$\mathcal{M}\ (0\leq\nu\leq m)$ and we denote a semi-Riemannian manifold by
$\mathcal{M}_{\nu}^{m}$. If $\nu=0$, then $\mathcal{M}$ is a Riemannian manifold.

The pair $(\nabla,g)$ is called a \textit{statistical structure} on
$\mathcal{M}$, if $\nabla$ is torsion-free and for vector fields $E,F,G$ on
$\mathcal{M}$%
\begin{equation}
(\nabla_{E}g)(F,G)=(\nabla_{F}g)(E,G) \label{1}%
\end{equation}
holds. (\ref{1}) is generally called \textit{Codazzi equation}. The triple
$(\mathcal{M},\nabla,g)$ is called a \textit{statistical manifold}. For the
statistical manifold $(\mathcal{M},\nabla,g)$, we define another affine
connection $\nabla^{\ast}$ by%
\begin{equation}
Eg(F,G)=g(\nabla_{E}F,G)+g(F,\nabla_{E}^{\ast}G). \label{2}%
\end{equation}
The affine connection $\nabla^{\ast}$ is called \textit{conjugate} or
\textit{dual} of $\nabla$ with respect to $g$. The affine connection
$\nabla^{\ast}$ is torsion-free and satisfies $(\nabla^{\ast})^{\ast}=\nabla$.
It is easy to see that $\hat{\nabla}=\frac{1}{2}(\nabla+\nabla^{\ast})$ is a
metric connection. The pair $(\nabla,g)$ is a statistical structure on
$\mathcal{M}$ if and only if so is $(\nabla^{\ast},g)$. Clearly, the triple
$(\mathcal{M},\nabla^{\ast},g)$ is statistical manifold.

We denote by $R$ and $R^{\ast}$ the curvature tensors on $\mathcal{M}$ with
respect to the affine connection $\nabla$ and its conjugate $\nabla^{\ast}$,
respectively. Then, we find
\begin{equation}
g(R(E,F)G,H)=-g(G,R^{\ast}(E,F)H), \label{2*}%
\end{equation}
where $R(E,F)G=[\nabla_{E},\nabla_{F}]G-\nabla_{\lbrack E,F]}G$. We put
\begin{equation}
S_{E}F=\nabla_{E}F-\nabla_{E}^{\ast}F. \label{2**}%
\end{equation}
Then $S_{E}F=S_{F}E$ and $g(S_{E}F,G)=g(F,S_{E}G)$ hold.

An almost complex structure on $\mathcal{M}$ is a tensor field $J$ of type
$(1,1)$ such that $J^{2}=-I$, where $I$ stands for the identity
transformation. An almost complex manifold is such a manifold with a fixed
almost complex structure. An almost complex manifold is necessarily orientable
and must have an even dimension. If $J$ preserves the metric $g$, that is,
\begin{equation}
g(JE,JF)=g(E,F), \label{10'}%
\end{equation}
then $(\mathcal{M},g,J)$ is an almost Hermitian manifold. Moreover, if $J$ is
parallel with respect to the Levi-Civita connection $\hat{\nabla}$, that is,
\begin{equation}
(\hat{\nabla}_{E}J)F=0, \label{10''}%
\end{equation}
then $(\mathcal{M},g,J)$ is called a K\"{a}hlerian manifold \cite{Yano}.

Let $(\mathcal{M},g,J)$ be a K\"{a}hlerian manifold and $\nabla$ an affine
connection of $\mathcal{M}$. We put $\omega(E,F)=g(E,JF)$ and $(\nabla
_{E}\omega)(F,G)=E\omega(F,G)-\omega(\nabla_{E}F,G)-\omega(F,\nabla_{E}G)$. If
$(\nabla,g)$ is a statistical structure and $\omega$ is a $\nabla$-parallel
2-form on $\mathcal{M}$, then $(\mathcal{M},\nabla,g,J)$ is called a
\textit{holomorphic statistical manifold }\cite{Mirjana}.

It is known that the following result \cite{Furuhata1}:

\noindent\textbf{Lemma A.} \textit{The following hold for a holomorphic
statistical manifold} $(\mathcal{M},\nabla,g,J)$:
\begin{align}
&  \nabla_{E}(JF)=J\nabla_{E}^{\ast}F,\label{18}\\
&  R(E,F)JG=JR^{\ast}(E,F)G. \label{18'}%
\end{align}
From (\ref{18}), we find $S_{E}(JF)=-J(S_{E}F).$

\begin{example}
\label{EX1}Let $%
%TCIMACRO{\U{211d} }%
%BeginExpansion
\mathbb{R}
%EndExpansion
_{2}^{4}$ be a smooth manifold with local coordinate system $(x_{1}%
,x_{2},x_{3},x_{4}),$ which admits the following almost complex structure $J:$%
\[
J=\left(
\begin{array}
[c]{cccc}%
0 & 0 & 1 & 0\\
0 & 0 & 0 & 1\\
-1 & 0 & 0 & 0\\
0 & -1 & 0 & 0
\end{array}
\right)  .
\]
The triple $(%
%TCIMACRO{\U{211d} }%
%BeginExpansion
\mathbb{R}
%EndExpansion
_{2}^{4},g,J)$ is an almost Hermitian manifold with
\[
g=\left(
\begin{array}
[c]{cccc}%
-1 & 0 & 0 & 0\\
0 & e^{-x_{2}} & 0 & 0\\
0 & 0 & -1 & 0\\
0 & 0 & 0 & e^{-x_{2}}%
\end{array}
\right)  .
\]
We set%
\begin{align*}
\nabla_{\partial_{1}}\partial_{1}  &  =-\nabla_{\partial_{3}}\partial
_{3}=-\partial_{2},\\
\nabla_{\partial_{1}}\partial_{2}  &  =\nabla_{\partial_{2}}\partial
_{1}=-\nabla_{\partial_{3}}\partial_{4}=-\nabla_{\partial_{4}}\partial
_{3}=e^{-x_{2}}\partial_{1}+e^{x_{1}}\partial_{4},\\
\nabla_{\partial_{1}}\partial_{3}  &  =\nabla_{\partial_{3}}\partial
_{1}=\partial_{4},\\
\nabla_{\partial_{1}}\partial_{4}  &  =\nabla_{\partial_{4}}\partial
_{1}=\nabla_{\partial_{2}}\partial_{3}=\nabla_{\partial_{3}}\partial
_{2}=e^{x_{1}}\partial_{2}-e^{-x_{2}}\partial_{3},\\
\nabla_{\partial_{2}}\partial_{2}  &  =-\nabla_{\partial_{4}}\partial
_{4}=-e^{x_{1}-x_{2}}\partial_{3},\\
\nabla_{\partial_{2}}\partial_{4}  &  =\nabla_{\partial_{4}}\partial
_{2}=-e^{x_{1}-x_{2}}\partial_{1}-\partial_{4,}%
\end{align*}
where $\partial_{i}=\partial/\partial x_{i}$ $(i=1,2,3,4).$ Then $(%
%TCIMACRO{\U{211d} }%
%BeginExpansion
\mathbb{R}
%EndExpansion
_{2}^{4},\nabla,g,J)$ is a holomorphic statistical manifold.
\end{example}

\section{\textbf{Holomorphic Statistical Submersions}}

Let $\mathcal{M}$ and $\mathcal{B}$ be semi-Riemannian manifolds. A surjective
mapping $\pi:\mathcal{M}\rightarrow\mathcal{B}$ is called a semi-Riemannian
submersion if $\pi$ has maximal rank and $\pi_{\ast}$ preserves lenghts of
horizontal vectors.

Let $\pi:\mathcal{M}\rightarrow\mathcal{B}$ be a semi-Riemannian submersion.
We put $\dim\mathcal{M}=m$ and $\dim\mathcal{B}=n$. For each point
$x\in\mathcal{B}$, semi-Riemannian submanifold $\pi^{-1}(x)$ with the induced
metric $\overline{g}$ is called a fiber and denoted by $\overline{\mathcal{M}%
}_{x}$ or $\overline{\mathcal{M}}$ simply. We notice that the dimension of
each fiber is always $m-n\,(=s)$. A vector field on $M$ is vertical if it is
always tangent to fibers, horizontal if always orthogonal to fibers. We denote
the vertical and horizontal subspace in the tangent space $T_{p}\mathcal{M}$
of the total space $\mathcal{M}$ by $\mathcal{V}_{p}(\mathcal{M})$ and
$\mathcal{H}_{p}(\mathcal{M})$ for each point $p\in\mathcal{M}$, and the
vertical and horizontal distributions in the tangent bundle $T\mathcal{M}$ of
$\mathcal{M}$ by $\mathcal{V}(\mathcal{M})$ and $\mathcal{H}(\mathcal{M})$,
respectively. Then $T\mathcal{M}$ is the direct sum of $\mathcal{V}%
(\mathcal{M})$ and $\mathcal{H}(\mathcal{M})$. The projection mappings are
denoted $\mathcal{V}:T\mathcal{M}\rightarrow\mathcal{V}(\mathcal{M})$ and
$\mathcal{H}:T\mathcal{M}\rightarrow\mathcal{H}(\mathcal{M})$ respectively. We
call a vector field $X$ on $\mathcal{M}$ projectable if there exists a vector
field $X_{\ast}$ on $\mathcal{B}$ such that $\pi_{\ast}(X_{p})=X_{\ast\pi(p)}$
for each $p\in\mathcal{M}$, and say that $X$ and $X_{\ast}$ are $\pi$-related.
Also, a vector field $X$ on $\mathcal{M}$ is called basic if it is projectable
and horizontal. Then, we have (\cite{Oneill}, \cite{Oneill 2}) \newline

\noindent\textbf{Lemma B.} \textit{If $X$ and $Y$ are basic vector fields on
}$\mathcal{M}$\textit{ which are }$\pi$\textit{-related to $X_{\ast}$ and
$Y_{\ast}$ on }$\mathcal{B}$\textit{, then}

\textit{i)}\ \ \textit{$g(X,Y)=\widetilde{g}(X_{\ast},Y_{\ast})\circ$}$\pi
$\textit{, where $g$ is the metric on $M$ and }$\widetilde{g}$\textit{ the
metric on }$\mathcal{B}$,

\textit{ii)}\ \ \textit{$\mathcal{H}[X,Y]$ is basic and is }$\pi
$\textit{-related to $[X_{\ast},Y_{\ast}]$.}

\vspace{.3cm}

Let $(\mathcal{M},\nabla,g)$ be a statistical manifold and $\pi:\mathcal{M}%
\rightarrow\mathcal{B}$ be a semi-Riemannian submersion. We denote the affine
connections of $\overline{\mathcal{M}}$ by $\overline{\nabla}$ and
$\overline{\nabla}^{\ast}$. Notice that $\overline{\nabla}_{U}V$ and
$\overline{\nabla}_{U}^{\ast}V$ are well-defined vertical vector fields on
$\mathcal{M}$ for vertical vector fields $U$ and $V$ on $\mathcal{M}$, more
precisely $\overline{\nabla}_{U}V=\mathcal{V}\nabla_{U}V$ and $\overline
{\nabla}_{U}^{\ast}V=\mathcal{V}\nabla_{U}^{\ast}V$. Moreover, $\overline
{\nabla}$ and $\overline{\nabla}^{\ast}$ are torsion-free and conjugate to
each other with respect to $\overline{g}$. Let $(\mathcal{M},\nabla,g)$ and
$\pi:\mathcal{M}\rightarrow\mathcal{B}$ be a statistical manifold and a
semi-Riemannian submersion, respectively. We call that $\pi:(\mathcal{M}%
,\nabla,g)\rightarrow(\mathcal{B},\widetilde{\nabla},\widetilde{g})$ is a
\textit{statistical submersion} if $\pi:\mathcal{M}\rightarrow\mathcal{B}$
satisfies $\pi_{\ast}(\nabla_{X}Y)_{p}=(\widetilde{\nabla}_{X_{\ast}}Y_{\ast
})_{\pi(p)}$ for basic vector fileds $X,Y$ and $p\in\mathcal{M}$. The letters
$U,V,W$ will always denote vertical vector fields, and $X,Y,Z$ horizontal
vector fields. The tensor fields $T$ and $A$ of type (1,2) defined by
\[
T_{E}F=\mathcal{H}\nabla_{\mathcal{V}E}\mathcal{V}F+\mathcal{V}\nabla
_{\mathcal{V}E}\mathcal{H}F,\hspace{1.5cm}A_{E}F=\mathcal{H}\nabla
_{\mathcal{H}E}\mathcal{V}F+\mathcal{V}\nabla_{\mathcal{H}E}\mathcal{H}F
\]
for vector fields $E$ and $F$ on $M$. Changing $\nabla$ to $\nabla^{\ast}$ in
the above equations, we set $T^{\ast}$ and $A^{\ast}$, respectively. Then we
find $(T^{\ast})^{\ast}=T$ and $(A^{\ast})^{\ast}=A$. For vertical vector
fields, $T$ and $T^{\ast}$ have the symmetry property. For $X,Y\in
\mathcal{H}(\mathcal{M})$ and $U,V\in\mathcal{V}(\mathcal{M})$, we obtain
\begin{equation}
g(T_{U}V,X)=-g(V,T_{U}^{\ast}X),\hspace{1.5cm}g(A_{X}Y,U)=-g(Y,A_{X}^{\ast}U).
\label{BB}%
\end{equation}
Thus, $T$ (resp. $A$) vanishes identically if and only if $T^{\ast}$ (resp.
$A^{\ast}$) vanishes identically. Since $A$ is related to the integrability of
$\mathcal{H}(\mathcal{M})$, if it is identically zero, then $\mathcal{H}%
(\mathcal{M})$ is integrable with respect to $\nabla$. Moreover, if $A$ and
$T$ vanish identically, then the total space is a locally product space of the
base space and the fiber. It is known that (\cite{Abe}) \newline

\noindent\textbf{Theorem C.} \textit{Let }$\pi:\mathcal{M}\rightarrow
\mathcal{B}$ \textit{ be a semi-Riemannian submersion. Then $(\mathcal{M}%
,\nabla,g)$ is a statistical manifold if and only if the following conditions
hold}:

\textit{i)}\ \ $\mathcal{H}S_{V}X=A_{X}V-A_{X}^{\ast}V$,

\textit{ii)}\ \ $\mathcal{V}S_{X}V=T_{V}X-T_{V}^{\ast}X$,

\textit{iii)}\ \ \textit{$(\overline{\mathcal{M}},\overline{\nabla}%
,\overline{g})$ is a statistical manifold for each $x\in$}$\mathcal{B}$,

\textit{iv)}\ \ \textit{$(\mathcal{B},\widetilde{\nabla},\widetilde{g})$ is a
statistical manifold}.

\vspace{.2cm}

For the statistical submersion $\pi:(\mathcal{M},\nabla,g)\rightarrow
(\mathcal{B},\widetilde{\nabla},\widetilde{g})$, we have the following Lemmas
(\cite{takano}) \newline

\noindent\textbf{Lemma D.} \textit{If $X$ and $Y$ are horizontal vector
fields, then $A_{X}Y=-A_{Y}^{\ast}X$.} \newline

\noindent\textbf{Lemma E.} \label{LEM1}\textit{For $X,Y\in\mathcal{H}%
(\mathcal{M})$ and $U,V\in\mathcal{V}(\mathcal{M})$ we have}
\[%
\begin{array}
[c]{ll}%
\nabla_{U}V=T_{U}V+\overline{\nabla}_{U}V, & \hspace{1.5cm}\nabla_{U}^{\ast
}V=T_{U}^{\ast}V+\overline{\nabla}_{U}^{\ast}V,\\[0.15cm]%
\nabla_{U}X=\mathcal{H}\nabla_{U}X+T_{U}X, & \hspace{1.5cm}\nabla_{U}^{\ast
}X=\mathcal{H}\nabla_{U}^{\ast}X+T_{U}^{\ast}X,\\[0.15cm]%
\nabla_{X}U=A_{X}U+\mathcal{V}\nabla_{X}U, & \hspace{1.5cm}\nabla_{X}^{\ast
}U=A_{X}^{\ast}U+\mathcal{V}\nabla_{X}^{\ast}U,\\[0.15cm]%
\nabla_{X}Y=\mathcal{H}\nabla_{X}Y+A_{X}Y, & \hspace{1.5cm}\nabla_{X}^{\ast
}Y=\mathcal{H}\nabla_{X}^{\ast}Y+A_{X}^{\ast}Y.
\end{array}
\]
\textit{Furthermore, if $X$ is basic, then $\mathcal{H}\nabla_{U}X=A_{X}U$ and
$\mathcal{H}\nabla_{U}^{\ast}X=A_{X}^{\ast}U$.}

\vspace{0.1cm} We define the covariant derivatives $\nabla T$ and $\nabla A$
by
\begin{align}
&  (\nabla_{E}T)_{F}G=\nabla_{E}(T_{F}G)-T_{\nabla_{E}F}G-T_{F}(\nabla
_{E}G),\nonumber\\
&  (\nabla_{E}A)_{F}G=\nabla_{E}(A_{F}G)-A_{\nabla_{E}F}G-A_{F}(\nabla
_{E}G)\nonumber
\end{align}
for $E,F,G\in T\mathcal{M}$. We change $\nabla$ to $\nabla^{\ast}$, then the
covariant derivatives $\nabla^{\ast}T,$ $\nabla^{\ast}A$ are defined simiraly.
We consider the curvature tensor on the statistical submersion. Let
$\overline{R}$ (resp. $\overline{R}^{\ast}$) be the curvature tensor with
respect to the induced affine connection $\overline{\nabla}$ (resp.
$\overline{\nabla}^{\ast}$ ) of each fiber. Also, let $\widetilde{R}%
(X,Y)Z\ $(resp. $\widetilde{R}^{\ast}(X,Y)Z$) be horizontal vector field such
that $\pi_{\ast}(\widetilde{R}(X,Y)Z)=\widetilde{R}(\pi_{\ast}X,\pi_{\ast
}Y)\pi_{\ast}Z$ (resp. $\pi_{\ast}(\widetilde{R}^{\ast}(X,Y)Z)=\widetilde
{R}^{\ast}(\pi_{\ast}X,\pi_{\ast}Y)\pi_{\ast}Z$) at each $p\in\mathcal{M}$,
where $\widetilde{R}$ (resp. $\widetilde{R}^{\ast}$) is the curvature tensor
on $\mathcal{B}$ of the affine connection $\widetilde{\nabla}$ (resp.
$\widetilde{\nabla}^{\ast}$). Then we have (\cite{takano}) \newline

\noindent\textbf{Theorem F.} \textit{If }$\pi$\textit{$:(\mathcal{M}%
,\nabla,g)\rightarrow(\mathcal{B},\widetilde{\nabla},\widetilde{g})$ is a
statistical submersion, then we get for $X,Y,Z,Z^{\prime}\in\mathcal{H}%
(\mathcal{M})$ and $U,V,W,W^{\prime}\in\mathcal{V}(\mathcal{M})$}
\begin{align}
&  g(R(U,V)W,W^{\prime})=g(\overline{R}(U,V)W,W^{\prime})+g(T_{U}W,T_{V}%
^{\ast}W^{\prime})-g(T_{V}W,T_{U}^{\ast}W^{\prime}),\nonumber\\
&  g(R(U,V)W,X)=g((\nabla_{U}T)_{V}W,X)-g((\nabla_{V}T)_{U}W,X),\nonumber\\
&  g(R(U,V)X,W)=g((\nabla_{U}T)_{V}X,W)-g((\nabla_{V}T)_{U}X,W),\nonumber\\
&  g(R(U,V)X,Y)=g((\nabla_{U}A)_{X}V,Y)-g((\nabla_{V}A)_{X}U,Y)+g(T_{U}%
X,T_{V}^{\ast}Y)-g(T_{V}X,T_{U}^{\ast}Y)\nonumber\\
&  \hspace{2.85cm}-g(A_{X}U,A_{Y}^{\ast}V)+g(A_{X}V,A_{Y}^{\ast}U),\nonumber\\
&  g(R(X,U)V,W)=g([\mathcal{V}\nabla_{X},\overline{\nabla}_{U}]V,W)-g(\nabla
_{\lbrack X,U]}V,W)-g(T_{U}V,A_{X}^{\ast}W)+g(T_{U}^{\ast}W,A_{X}%
V),\nonumber\\
&  g(R(X,U)V,Y)=g((\nabla_{X}T)_{U}V,Y)-g((\nabla_{U}A)_{X}V,Y)+g(A_{X}%
U,A_{Y}^{\ast}V)-g(T_{U}X,T_{V}^{\ast}Y),\nonumber\\
&  g(R(X,U)Y,V)=g((\nabla_{X}T)_{U}Y,V)-g((\nabla_{U}A)_{X}Y,V)+g(T_{U}%
X,T_{V}Y)-g(A_{X}U,A_{Y}V),\nonumber\\
&  g(R(X,U)Y,Z)=g((\nabla_{X}A)_{Y}U,Z)-g(T_{U}X,A_{Y}^{\ast}Z)-g(T_{U}%
Y,A_{X}^{\ast}Z)+g(A_{X}Y,T_{U}^{\ast}Z),\nonumber\\
&  g(R(X,Y)U,V)=g([\mathcal{V}\nabla_{X},\mathcal{V}\nabla_{Y}]U,V)-g(\nabla
_{\lbrack X,Y]}U,V)+g(A_{X}U,A_{Y}^{\ast}V)-g(A_{Y}U,A_{X}^{\ast
}V),\nonumber\\
&  g(R(X,Y)U,Z)=g((\nabla_{X}A)_{Y}U,Z)-g((\nabla_{Y}A)_{X}U,Z)+g(T_{U}^{\ast
}Z,\theta_{X}Y),\nonumber\\
&  g(R(X,Y)Z,U)=g((\nabla_{X}A)_{Y}Z,U)-g((\nabla_{Y}A)_{X}Z,U)-g(T_{U}%
Z,\theta_{X}Y),\nonumber\\
&  g(R(X,Y)Z,Z^{\prime})=g(\widehat{R}(X,Y)Z,Z^{\prime})-g(A_{Y}Z,A_{X}^{\ast
}Z^{\prime})+g(A_{X}Z,A_{Y}^{\ast}Z^{\prime})+g(\theta_{X}Y,A_{Z}^{\ast
}Z^{\prime}),\nonumber
\end{align}
\textit{where we put $\theta_{X}=A_{X}+A_{X}^{\ast}$.} \newline

\noindent\textbf{Remark G.} \vspace{0.2cm}We find $\mathcal{V[}X,Y]=$%
\textit{$\theta_{X}Y.$}

Let $(\mathcal{M},\nabla,g,J)$ be a holomorphic statistical manifold and
$(\mathcal{B},\widetilde{\nabla},\widetilde{g})$ be a statistical manifold.
The statistical submersion $\pi:(\mathcal{M},\nabla,g,J)\rightarrow
(\mathcal{B},\widetilde{\nabla},\widetilde{g})$ is called a
\textit{holomorphic statistical submersion}. For $X\in\mathcal{H}%
(\mathcal{M})$ and $U\in\mathcal{V}(\mathcal{M})$ we put
\begin{equation}
JX=PX+FX,\hspace{2cm}JU=tU+fU, \label{A}%
\end{equation}
where $PX,tU\in\mathcal{H}(\mathcal{M})$ and $FX,fU\in\mathcal{V}%
(\mathcal{M}).$

From $J^{2}=-I$, we get
\[
P^{2}=-I-tF,\hspace{1cm}FP+fF=0,\hspace{1cm}Pt+tf=0,\hspace{1cm}f^{2}=-I-Ft.
\]
Because of $g(JE,G)+g(E,JG)=0$ for $E,G\in T\mathcal{M}$, we find
\begin{align}
&  g(PY,Z)+g(Y,PZ)=0,\label{AA}\\
&  g(FX,U)+g(X,tU)=0,\label{AB}\\
&  g(fV,W)+g(V,fW)=0. \label{AC}%
\end{align}
Moreover, we obtain
\begin{align}
&  g((\mathcal{H}\nabla_{X}P)Y,Z)+g(Y,(\mathcal{H}\nabla_{X}^{\ast
}P)Z)=0,\nonumber\\
&  g((\mathcal{H}\nabla_{U}P)Y,Z)+g(Y,(\mathcal{H}\nabla_{U}^{\ast
}P)Z)=0,\nonumber\\
&  g((\mathcal{V}\nabla_{X}f)V,W)+g(V,(\mathcal{V}\nabla_{X}^{\ast
}f)W)=0,\nonumber\\
&  g((\overline{\nabla}_{U}f)V,W)+g(V,(\overline{\nabla}_{U}^{\ast
}f)W)=0.\nonumber
\end{align}
Hence we have \newline

\noindent\textbf{Lemma 3.1.} \ \label{L1}\ \textit{If }$\pi:(\mathcal{M}%
,\nabla,g,J)\rightarrow(\mathcal{B},\widetilde{\nabla},\widetilde{g})$\textit{
is a holomorphic statistical submersion, then we have}

\textit{i)}\ \ \textit{$\mathcal{H}\nabla_{X}P=0$ $($resp. $\mathcal{H}%
\nabla_{U}P=0)$ is equivalent to $\mathcal{H}\nabla_{X}^{\ast}P=0$ $($resp.
$\mathcal{H}\nabla_{U}^{\ast}P=0)$}.

\textit{ii)}\ \ \textit{$\mathcal{V}\nabla_{X}f=0$ $($resp. $\overline{\nabla
}_{U}f=0)$ is equivalent to $\mathcal{V}\nabla_{X}^{\ast}f=0$ $($resp.
$\overline{\nabla}_{U}^{\ast}f=0)$}.

\vspace{.2cm}

Using (\ref{18}), we can get \newline

\noindent\textbf{Lemma 3.2.} \label{L2}\textit{Let }$\pi:(\mathcal{M}%
,\nabla,g,J)\rightarrow(\mathcal{B},\widetilde{\nabla},\widetilde{g})$\textit{
be a holomorphic statistical submersion. Then we have}%
\begin{align}
&  \mathcal{H}\nabla_{U}(tV)+T_{U}(fV)=P(T_{U}^{\ast}V)+t(\overline{\nabla
}_{U}^{\ast}V),\label{K26}\\
&  T_{U}(tV)+\overline{\nabla}_{U}(fV)=F(T_{U}^{\ast}V)+f(\overline{\nabla
}_{U}^{\ast}V),\label{K27}\\
&  \mathcal{H}\nabla_{U}(PX)+T_{U}(FX)=P(\mathcal{H}\nabla_{U}^{\ast
}X)+t(T_{U}^{\ast}X),\label{K28}\\
&  T_{U}(PX)+\overline{\nabla}_{U}(FX)=F(\mathcal{H}\nabla_{U}^{\ast
}X)+f(T_{U}^{\ast}X),\label{K29}\\
&  \mathcal{H}\nabla_{X}(tU)+A_{X}(fU)=P(A_{X}^{\ast}U)+t(\mathcal{V}%
\nabla_{X}^{\ast}U),\label{K30}\\
&  A_{X}(tU)+\mathcal{V}\nabla_{X}(fU)=F(A_{X}^{\ast}U)+f(\mathcal{V}%
\nabla_{X}^{\ast}U),\label{K31}\\
&  \mathcal{H}\nabla_{X}(PY)+A_{X}(FY)=P(\mathcal{H}\nabla_{X}^{\ast
}Y)+t(A_{X}^{\ast}Y),\label{K32}\\
&  A_{X}(PY)+\mathcal{V}\nabla_{X}(FY)=F(\mathcal{H}\nabla_{X}^{\ast
}Y)+f(A_{X}^{\ast}Y). \label{K33}%
\end{align}
\textit{Furthermore, if $X$ is basic, then $\mathcal{H}\nabla_{U}^{\ast
}X=A_{X}^{\ast}U$.} \newline

\noindent\textbf{Corollary 3.3.} \textit{Let }$\pi:(\mathcal{M},\nabla
,g,J)\rightarrow(\mathcal{B},\widetilde{\nabla},\widetilde{g})$\textit{ a
holomorphic statistical submersion. Then we get}%
\begin{align}
&  T_{U}^{\ast}V=-P\{\mathcal{H}\nabla_{U}(tV)+T_{U}(fV)\}-t\{T_{U}%
(tV)+\overline{\nabla}_{U}(fV)\},\label{K34}\\
&  \overline{\nabla}_{U}^{\ast}V=-F\{\mathcal{H}\nabla_{U}(tV)+T_{U}%
(fV)\}-f\{T_{U}(tV)+\overline{\nabla}_{U}(fV)\},\label{K35}\\
&  \mathcal{H}\nabla_{U}^{\ast}X=-P\{\mathcal{H}\nabla_{U}(PX)+T_{U}%
(FX)\}-t\{T_{U}(PX)+\overline{\nabla}_{U}(FX)\},\label{K36}\\
&  T_{U}^{\ast}X=-F\{\mathcal{H}\nabla_{U}(PX)+T_{U}(FX)\}-f\{T_{U}%
(PX)+\overline{\nabla}_{U}(FX)\},\label{K37}\\
&  A_{X}^{\ast}U=-P\{\mathcal{H}\nabla_{X}(tU)+A_{X}(fU)\}-t\{A_{X}%
(tU)+\mathcal{V}\nabla_{X}(fU)\},\label{K38}\\
&  \mathcal{V}\nabla_{X}^{\ast}U=-F\{\mathcal{H}\nabla_{X}(tU)+A_{X}%
(fU)\}-f\{A_{X}(tU)+\mathcal{V}\nabla_{X}(fU)\},\label{K39}\\
&  \mathcal{H}\nabla_{X}^{\ast}Y=-P\{\mathcal{H}\nabla_{X}(PY)+A_{X}%
(FY)\}-t\{A_{X}(PY)+\mathcal{V}\nabla_{X}(FY)\},\label{K40}\\
&  A_{X}^{\ast}Y=-F\{\mathcal{H}\nabla_{X}(PY)+A_{X}(FY)\}-f\{A_{X}%
(PY)+\mathcal{V}\nabla_{X}(FY)\}. \label{K41}%
\end{align}
We put%
\begin{align}
&  (\overline{\nabla}_{U}f)V=\overline{\nabla}_{U}(fV)-f(\overline{\nabla}%
_{U}V),\nonumber\\
&  (\mathcal{H}\nabla_{U}P)X=\mathcal{H}\nabla_{U}(PX)-P(\mathcal{H}\nabla
_{U}X),\nonumber\\
&  (\mathcal{V}\nabla_{X}f)U=\mathcal{V}\nabla_{X}(fU)-f(\mathcal{V}\nabla
_{X}U),\nonumber\\
&  (\mathcal{H}\nabla_{X}P)Y=\mathcal{H}\nabla_{X}(PY)-P(\mathcal{H}\nabla
_{X}Y).\nonumber
\end{align}
From (\ref{K27}), (\ref{K28}), (\ref{K31}) and (\ref{K32}), we obtain \newline

\noindent\textbf{Corollary 3.4.} \textit{Let }$\pi:(\mathcal{M},\nabla
,g,J)\rightarrow(\mathcal{B},\widetilde{\nabla},\widetilde{g})$\textit{ a
holomorphic statistical submersion. Then we get}%
\begin{align*}
&  (\overline{\nabla}_{U}f)V=-f(\mathcal{V}(S_{U}V))+F(T_{U}^{\ast}%
V)-T_{U}(tV),\\
&  (\mathcal{H}\nabla_{U}P)X=-P(\mathcal{H}(S_{U}X))+t(T_{U}^{\ast}%
X)-T_{U}(FX),\\
&  (\mathcal{V}\nabla_{X}f)U=-f(\mathcal{V}(S_{X}U))+F(A_{X}^{\ast}%
U)-A_{X}(tU),\\
&  (\mathcal{H}\nabla_{X}P)Y=-P(\mathcal{H}(S_{X}Y))+t(A_{X}^{\ast}%
Y)-A_{X}(FY).
\end{align*}
\newline

\noindent\textbf{Corollary 3.5.} \textit{Let }$\pi:(\mathcal{M},\nabla
,g,J)\rightarrow(\mathcal{B},\widetilde{\nabla},\widetilde{g})$\textit{ a
holomorphic statistical submersion. Then we get}%
\begin{align*}
&  T_{U}(tV)+f(\overline{\nabla}_{U}V)=F(T_{U}^{\ast}V)+f(\overline{\nabla
}_{U}^{\ast}V)\text{ \ \ \ \ \ \ } & if  &  \text{ \ }\overline{\nabla}%
_{U}f=0,\\
&  P(\mathcal{H}\nabla_{U}X)+T_{U}(FX)=P(\mathcal{H}\nabla_{U}^{\ast
}X)+t(T_{U}^{\ast}X)\text{ \ \ \ \ } & if  &  \text{ \ }\mathcal{H}\nabla
_{U}P=0,\\
&  A_{X}(tU)+f(\mathcal{V}\nabla_{X}U)=F(A_{X}^{\ast}U)+f(\mathcal{V}%
\nabla_{X}^{\ast}U)\text{ \ \ \ } & if  &  \text{ \ }\mathcal{V}\nabla
_{X}f=0,\\
&  P(\mathcal{H}\nabla_{X}Y)+A_{X}(FY)=P(\mathcal{H}\nabla_{X}^{\ast
}Y)+t(A_{X}^{\ast}Y)\text{ \ \ \ } & if  &  \text{ }\mathcal{H}\nabla_{X}P=0.
\end{align*}
Now, we can give an example of the holomorphic statistical submersion:

\begin{example}
\label{EX2}Let $(%
%TCIMACRO{\U{211d} }%
%BeginExpansion
\mathbb{R}
%EndExpansion
_{1}^{2},\widetilde{g})$ be a semi-Riemannian manifold with local coordinate
system $(x_{1},x_{2}),$ where $\widetilde{g}=\left(
\begin{array}
[c]{cc}%
-1 & 0\\
0 & e^{-x_{2}}%
\end{array}
\right)  .$ If we put%
\[
\widetilde{\nabla}_{\partial_{1\ast}}\partial_{1\ast}=-\partial_{2\ast},\text{
\ \ \ \ \ }\widetilde{\nabla}_{\partial_{1\ast}}\partial_{2\ast}%
=\widetilde{\nabla}_{\partial_{2\ast}}\partial_{1\ast}=e^{-x_{2}}%
\partial_{1\ast},\text{ \ \ \ }\widetilde{\nabla}_{\partial_{2\ast}}%
\partial_{2\ast}=0,
\]
then $(%
%TCIMACRO{\U{211d} }%
%BeginExpansion
\mathbb{R}
%EndExpansion
_{1}^{2},\widetilde{\nabla},\widetilde{g})$ is a statistical manifold, where
$\partial_{i\ast}=\partial/\partial x_{i}$ $(i=1,2).$ Considering the
holomorphic statistical manifold $(%
%TCIMACRO{\U{211d} }%
%BeginExpansion
\mathbb{R}
%EndExpansion
_{2}^{4},\nabla,g,J)$ given in Example\ref{EX1}, we define a holomorphic
statistical submersion $\pi:(%
%TCIMACRO{\U{211d} }%
%BeginExpansion
\mathbb{R}
%EndExpansion
_{2}^{4},\nabla,g,J)\rightarrow(%
%TCIMACRO{\U{211d} }%
%BeginExpansion
\mathbb{R}
%EndExpansion
_{1}^{2},\widetilde{\nabla},\widetilde{g})$ by%
\[
\pi(x_{1},x_{2},x_{3},x_{4})=(x_{1},x_{2}).
\]
Moreover, for $\partial_{1},\partial_{2}\in\Gamma(\mathcal{H})$ and
$\partial_{3},\partial_{4}\in\Gamma(\mathcal{V}),$ we get%
\[%
\begin{array}
[c]{ll}%
T_{\partial_{3}}\partial_{3}=\partial_{2}, & \overline{\nabla}_{\partial_{3}%
}\partial_{3}=0,\\
T_{\partial_{3}}\partial_{4}=T_{\partial_{4}}\partial_{3}=-e^{-x_{2}}%
\partial_{1}, & \overline{\nabla}_{\partial_{3}}\partial_{4}=\overline{\nabla
}_{\partial_{4}}\partial_{3}=-e^{x_{1}}\partial_{4},\\
T_{\partial_{4}}\partial_{4}=0, & \overline{\nabla}_{\partial_{4}}\partial
_{4}=e^{x_{1}-x_{2}}\partial_{3},\\
\mathcal{H}\nabla_{\partial_{3}}\partial_{1}=0, & T_{\partial_{3}}\partial
_{1}=\partial_{4},\\
\mathcal{H}\nabla_{\partial_{3}}\partial_{2}=\mathcal{H}\nabla_{\partial_{4}%
}\partial_{1}=e^{x_{1}}\partial_{2}, & T_{\partial_{3}}\partial_{2}%
=T_{\partial_{4}}\partial_{1}=-e^{-x_{2}}\partial_{3},\\
\mathcal{H}\nabla_{\partial_{4}}\partial_{2}=-e^{x_{1}-x_{2}}\partial_{1}, &
T_{\partial_{4}}\partial_{2}=-\partial_{4},\\
A_{\partial_{1}}\partial_{3}=0, & \mathcal{V}\nabla_{\partial_{1}}\partial
_{3}=\partial_{4},\\
A_{\partial_{1}}\partial_{4}=A_{\partial_{2}}\partial_{3}=e^{x_{1}}%
\partial_{2}, & \mathcal{V}\nabla_{\partial_{1}}\partial_{4}=\mathcal{V}%
\nabla_{\partial_{2}}\partial_{3}=-e^{-x_{2}}\partial_{3},\\
A_{\partial_{2}}\partial_{4}=-e^{x_{1}-x_{2}}\partial_{1}, & \mathcal{V}%
\nabla_{\partial_{2}}\partial_{4}=-\partial_{4},\\
\mathcal{H}\nabla_{\partial_{1}}\partial_{1}=-\partial_{2}, & A_{\partial_{1}%
}\partial_{1}=0,\\
\mathcal{H}\nabla_{\partial_{1}}\partial_{2}=\mathcal{H}\nabla_{\partial_{2}%
}\partial_{1}=e^{-x_{2}}\partial_{1}, & A_{\partial_{1}}\partial
_{2}=A_{\partial_{2}}\partial_{1}=e^{x_{1}}\partial_{4},\\
\mathcal{H}\nabla_{\partial_{2}}\partial_{2}=0, & A_{\partial_{2}}\partial
_{2}=-e^{x_{1}-x_{2}}\partial_{3}.
\end{array}
\]

\end{example}

Let $\pi:(\mathcal{M},\nabla,g,J)\rightarrow(\mathcal{B},\widetilde{\nabla
},\widetilde{g})$ be a holomorphic statistical submersion. We consider the
curvature with respect to the affine connection $\nabla$ of the total space
satisfies
\begin{equation}
R(E,F)G=\frac{\,c\,}{4}\{g(F,G)E-g(E,G)F+g(JF,G)JE-g(JE,G)JF+2g(E,JF)JG\}
\label{*}%
\end{equation}
for $E,\,F,\,G\in T\mathcal{M}$, where $c$ is a constant. From (\ref{18'}), we
find
\begin{align}
&  g(\overline{R}(U,V)W,W^{\prime})+g(T_{U}W,T_{V}^{\ast}W^{\prime}%
)-g(T_{V}W,T_{U}^{\ast}W^{\prime})\\
&  =\frac{\,c\,}{4}\{g(V,W)g(U,W^{\prime})-g(U,W)g(V,W^{\prime}%
)+g(fV,W)g(fU,W^{\prime})\nonumber\\
&  \hspace{0.8cm}-g(fU,W)g(fV,W^{\prime})+2g(U,fV)g(fW,W^{\prime
})\},\nonumber\\
&  g((\nabla_{U}T)_{V}W,X)-g((\nabla_{V}T)_{U}W,X)\label{R2}\\
&  =\frac{\,c\,}{4}%
\{g(fV,W)g(tU,X)-g(fU,W)g(tV,X)+2g(U,fV)g(tW,X)\},\nonumber\\
&  g((\nabla_{U}T)_{V}X,W)-g((\nabla_{V}T)_{U}X,W)\label{R22}\\
&  =\frac{\,c\,}{4}%
\{g(tV,X)g(fU,W)-g(tU,X)g(fV,W)+2g(U,fV)g(FX,W)\},\nonumber\\
&  g((\nabla_{U}A)_{X}V,Y)-g((\nabla_{V}A)_{X}U,Y)+g(T_{U}X,T_{V}^{\ast
}Y)-g(T_{V}X,T_{U}^{\ast}Y)\\
&  -g(A_{X}U,A_{Y}^{\ast}V)+g(A_{X}V,A_{Y}^{\ast}U)\nonumber\\
&  =\frac{\,c\,}{4}%
\{g(tV,X)g(tU,Y)-g(tU,X)g(tV,Y)+2g(U,fV)g(PX,Y)\},\nonumber\\
&  g([\mathcal{V}\nabla_{X},\overline{\nabla}_{U}]V,W)-g(\nabla_{\lbrack
X,U]}V,W)-g(T_{U}V,A_{X}^{\ast}W)+g(A_{X}V,T_{U}^{\ast}W)\nonumber\\
&  =\frac{\,c\,}{4}%
\{g(fU,V)g(FX,W)-g(FX,V)g(fU,W)+2g(X,tU)g(fV,W)\},\nonumber\\
&  g((\nabla_{X}T)_{U}V,Y)-g((\nabla_{U}A)_{X}V,Y)+g(A_{X}U,A_{Y}^{\ast
}V)-g(T_{U}X,T_{V}^{\ast}Y)\\
&  =\frac{\,c\,}{4}\{g(U,V)g(X,Y)+g(fU,V)g(PX,Y)-g(FX,V)g(tU,Y)\nonumber\\
&  \hspace{0.85cm}+2g(X,tU)g(tV,Y)\},\nonumber\\
&  g((\nabla_{X}T)_{U}Y,V)-g((\nabla_{U}A)_{X}Y,V)-g(A_{X}U,A_{Y}%
V)+g(T_{U}X,T_{V}Y)\\
&  =\frac{\,c\,}{4}\{-g(X,Y)g(U,V)+g(tU,Y)g(FX,V)-g(PX,Y)g(fU,V)\nonumber\\
&  \hspace{0.85cm}+2g(X,tU)g(FY,V)\},\nonumber\\
&  g((\nabla_{X}A)_{Y}U,Z)-g(T_{U}X,A_{Y}^{\ast}Z)-g(T_{U}Y,A_{X}^{\ast
}Z)+g(A_{X}Y,T_{U}^{\ast}Z)\label{R3}\\
&  =\frac{\,c\,}{4}%
\{g(tU,Y)g(PX,Z)-g(PX,Y)g(tU,Z)+2g(X,tU)g(PY,Z)\},\nonumber\\
&  g([\mathcal{V}\nabla_{X},\mathcal{V}\nabla_{Y}]U,V)-g(\nabla_{\lbrack
X,Y]}U,V)+g(A_{X}U,A_{Y}^{\ast}V)-g(A_{Y}U,A_{X}^{\ast}V)\\
&  =\frac{\,c\,}{4}%
\{g(FY,U)g(FX,V)-g(FX,U)g(FY,V)+2g(X,PY)g(fU,V)\},\nonumber\\
&  g((\nabla_{X}A)_{Y}U,Z)-g((\nabla_{Y}A)_{X}U,Z)+g(T_{U}^{\ast}Z,\theta
_{X}Y)\\
&  =\frac{\,c\,}{4}%
\{g(FY,U)g(PX,Z)-g(FX,U)g(PY,Z)+2g(X,PY)g(tU,Z)\},\nonumber\\
&  g((\nabla_{X}A)_{Y}Z,U)-g((\nabla_{Y}A)_{X}Z,U)-g(T_{U}Z,\theta_{X}Y)\\
&  =\frac{\,c\,}{4}%
\{g(PY,Z)g(FX,U)-g(PX,Z)g(FY,U)+2g(X,PY)g(FZ,U)\},\nonumber\\
&  g(\widehat{R}(X,Y)Z,Z^{\prime})-g(A_{Y}Z,A_{X}^{\ast}Z^{\prime}%
)+g(A_{X}Z,A_{Y}^{\ast}Z^{\prime})+g(\theta_{X}Y,A_{Z}^{\ast}Z^{\prime})\\
&  =\frac{\,c\,}{4}\{g(Y,Z)g(X,Z^{\prime})-g(X,Z)g(Y,Z^{\prime}%
)+g(PY,Z)g(PX,Z^{\prime})\nonumber\\
&  \hspace{0.85cm}-g(PX,Z)g(PY,Z^{\prime})+2g(X,PY)g(PZ,Z^{\prime})\}\nonumber
\end{align}
for $X,\,Y,\,Z,\,Z^{\prime}\in\mathcal{H}(M)$ and $U,\,V,\,W,\,W^{\prime}%
\in\mathcal{V}(M)$, where $\theta_{X}Y=A_{X}Y+A_{X}^{\ast}Y$.

\vspace{.2cm}

Let $\pi:(\mathcal{M},\nabla,g,J)\rightarrow(\mathcal{B},\widetilde{\nabla
},\widetilde{g})$ be a holomorphic statistical submersion with isometric
fiber, that is, $T=0$. We get from (\ref{R2})
\[
c\{g(fV,W)g(tU,X)-g(fU,W)g(tV,X)+2g(U,fV)g(tW,X)\}=0.
\]
Thus we find $c=0$ or
\begin{equation}
g(fV,W)g(tU,X)-g(fU,W)g(tV,X)+2g(U,fV)g(tW,X)=0. \label{S}%
\end{equation}
Because of (\ref{S}), we get $g(FX,fV)=0$ which yields that $fF=0$ and $tf=0$.
Changing $W$ (resp. $V$) to $fW$ (resp. $fV$), equation (\ref{S}) are
\begin{align}
&  g(f^{2}V,W)g(tU,X)-g(f^{2}U,W)g(tV,X)=0,\nonumber\\
&  g(f^{2}V,W)g(tU,X)+2g(U,f^{2}V)g(tW,X)=0,\nonumber
\end{align}
respectively. Thus we obtain $g(f^{2}U,W)g(tV,X)=0$ which means that $t=0$ or
$f^{2}=0$. For each $p\in\mathcal{M}$, we denote by $\left\{  U_{1}%
,...,U_{s}\right\}  $ local orthonomal bases of $\mathcal{V}_{p}(\mathcal{M}%
)$, where $s=$dim$\overline{\mathcal{M}}$. When $f^{2}=0,$ we get $\left\Vert
f\right\Vert ^{2}=\sum\varepsilon_{\alpha}g(fU_{\alpha},fU_{\alpha}%
)=-\sum\varepsilon_{\alpha}g(f^{2}U_{\alpha},U_{\alpha})=0,$ that is, $f=0$.
Hence we have \newline

\noindent\textbf{Theorem 3.6.} \vspace{0.2cm} \textit{Let }$\pi:(\mathcal{M}%
,\nabla,g,J)\rightarrow(\mathcal{B},\widetilde{\nabla},\widetilde{g})$\textit{
be a holomorphic statistical submersion with isometric fiber. If the total
space satisfies the condition $(\ref{*})$, then}

\textit{i)}\ \ \textit{the total space is flat, or}

\textit{ii)}\ \ \textit{each fiber is an invariant submanifold of
}$\mathcal{M}$ \textit{satisfying the condition $(\ref{*})$, or}

\textit{iii)}\ \ \textit{each fiber is an anti-invariant submanifold of
}$\mathcal{M}$\textit{ which is of constant curvature $\frac{\,c\,}{4}$.}

\section{\textbf{Anti-Invariant Holomorphic Statistical Submersions}}

The holomorphic statistical submersion $\pi:(\mathcal{M},\nabla
,g,J)\rightarrow(\mathcal{B},\widetilde{\nabla},\widetilde{g})$\textit{ } is
called an \textit{anti-invariant holomorphic statistical submersion} if
$\overline{\mathcal{M}}$ is an anti-invariant submanifold of $\mathcal{M}$,
namely, $J(\mathcal{V}(\mathcal{M}))\subset\mathcal{H}(\mathcal{M})$. Thus, we
find $P^{2}=-I-tF,\ FP=0,\ Pt=0$ and $Ft=-I$. We assume $\mathcal{H}\nabla
_{U}P=0$. Then we get from Corollary 3.4
\begin{equation}
P(\mathcal{H(}S_{U}X))-t(T_{U}^{\ast}X)+T_{U}(FX)=0. \label{AA1}%
\end{equation}
If we operate $F$ to (\ref{AA1}), we obtain
\begin{equation}
T_{U}^{\ast}X=-F(T_{U}(FX)). \label{AA2}%
\end{equation}
Hence we have from $(T^{\ast})^{\ast}=T$ and (\ref{BB})

\noindent\textbf{Lemma 4.1.} \label{L2*}\textit{Let }$\pi:(\mathcal{M}%
,\nabla,g,J)\rightarrow(\mathcal{B},\widetilde{\nabla},\widetilde{g})$
be\textit{ an anti-invariant holomorphic statistical submersion. If
$\mathcal{H}\nabla_{U}P=0$, then we get}
\[%
\begin{array}
[c]{ll}%
T_{U}X=-F(T_{U}^{\ast}(FX)), & \hspace{1cm}T_{U}^{\ast}X=-F(T_{U}%
(FX)),\\[0.15cm]%
T_{U}V=-t(T_{U}^{\ast}(tV)), & \hspace{1cm}T_{U}^{\ast}V=-t(T_{U}(tV)).
\end{array}
\]

\noindent\textbf{Corollary 4.2.} \label{CR}\textit{Let }$\pi:(\mathcal{M}%
,\nabla,g,J)\rightarrow(\mathcal{B},\widetilde{\nabla},\widetilde{g})$
be\textit{ an anti-invariant holomorphic statistical submersion. If
$\mathcal{H}\nabla_{U}P=0$, then we find}
\[%
\begin{array}
[c]{ll}%
t(T_{U}X)=T_{U}^{\ast}(FX), & \hspace{1cm}t(T_{U}^{\ast}X)=T_{U}(FX),\\[0.1cm]%
P(T_{U}V)=0, & \hspace{1cm}P(T_{U}^{\ast}V)=0,\\[0.1cm]%
F(T_{U}V)=T_{U}^{\ast}(tV)=T_{V}^{\ast}(tU), & \hspace{1cm}F(T_{U}^{\ast
}V)=T_{U}(tV)=T_{V}(tU),\\[0.1cm]%
T_{U}(PX)=0, & \hspace{1cm}T_{U}^{\ast}(PX)=0,\\[0.1cm]%
P(\mathcal{H}\nabla_{U}X)=P(\mathcal{H}\nabla_{U}^{\ast}X). &
\end{array}
\]

Using (\ref{2}) and (\ref{10''}), we get%
\begin{align*}
g(\nabla_{U}V,X)  &  =g(J(\nabla_{U}V),JX)=g(\nabla_{U}^{\ast}(JV),JX)\\
&  =Ug(V,X)-g(tV,\nabla_{U}(PX))-g(tV,\nabla_{U}(FX))\\
&  =g(V,F((\mathcal{H}\nabla_{U}P)X))+g(V,F(T_{U}(FX))).
\end{align*}

Then we find $F((\mathcal{H}\nabla_{U}P)X)+T_{U}^{\ast}X+F(T_{U}(FX))=0.$ Thus
we find

\noindent\textbf{Lemma 4.3.} \label{L3}\textit{Let }$\pi:(\mathcal{M}%
,\nabla,g,J)\rightarrow(\mathcal{B},\widetilde{\nabla},\widetilde{g})$
be\textit{ an anti-invariant holomorphic statistical submersion. If
}$\mathcal{V}$\textit{ is a totally geodesic foliation on }$\mathcal{M}%
$\textit{, then} \textit{$\mathcal{H}\nabla_{X}P=0$ holds.}

Next, if \textit{$\mathcal{H}\nabla_{X}P=0$}, then we get from Corollary 3.4
\begin{equation}
P(\mathcal{H(}S_{X}Y\mathcal{)})-t(A_{X}^{\ast}Y)+A_{X}(FY)=0. \label{aaa'}%
\end{equation}

Operating $F$ to (\ref{aaa'})%
\begin{equation}
A_{X}^{\ast}Y=F(A_{X}(FY)). \label{aaa''}%
\end{equation}
Hence we have from $(A^{\ast})^{\ast}=A$ and (\ref{BB})

\noindent\textbf{Lemma 4.4.} \label{C1}{\textit{$\pi:(\mathcal{M}%
,\nabla,g,J)\rightarrow(\mathcal{B},\widetilde{\nabla},\widetilde{g})$ be an
anti-invariant holomorphic statistical submersion. If $\mathcal{H}\nabla
_{X}P=0$, then we get}}
\[%
\begin{array}
[c]{ll}%
A_{X}Y=-F(A_{X}^{\ast}(FY)), & \hspace{1.5cm}A_{X}^{\ast}Y=-F(A_{X}%
(FY)),\\[0.2cm]%
A_{X}U=-t(A_{X}^{\ast}(tU)), & \hspace{1.5cm}A_{X}^{\ast}U=-t(A_{X}(tU)).
\end{array}
\]

\noindent\textbf{Corollary 4.5.} \textit{Let }$\pi:(\mathcal{M},\nabla
,g,J)\rightarrow(\mathcal{B},\widetilde{\nabla},\widetilde{g})$\textit{ be an
anti-invariant holomorphic statistical submersion. If $\mathcal{H}\nabla
_{X}P=0,$ then we find}
\[%
\begin{array}
[c]{ll}%
\mathit{F(A_{X}U)=A_{X}^{\ast}(tU)}, & \hspace{1.5cm}\mathit{F(A_{X}^{\ast
}U)=A_{X}(tU)},\\[0.2cm]%
A_{X}(PY)=0, & \hspace{1.5cm}A_{X}^{\ast}(PY)=0,\\[0.2cm]%
P(A_{X}U)=0, & \hspace{1.5cm}P(A_{X}^{\ast}U)=0,\\[0.2cm]%
t(A_{X}Y)=A_{X}^{\ast}(FY), & \hspace{1.5cm}t(A_{X}^{\ast}Y)=A_{X}%
(FY),\\[0.2cm]%
A_{PY}U=0, & \hspace{1.5cm}A_{PY}^{\ast}U=0,\\[0.2cm]%
P(\mathcal{H}\nabla_{X}Y)=P(\mathcal{H}\nabla_{X}^{\ast}Y). &
\end{array}
\]
Owing to (\ref{K35}), (\ref{K39}) and $(\nabla^{\ast})^{\ast}=\nabla$, we get

\noindent\textbf{Lemma 4.6.} \textit{Let }$\pi:(\mathcal{M},\nabla
,g,J)\rightarrow(\mathcal{B},\widetilde{\nabla},\widetilde{g}$\textit{$)$ be
an anti-invariant holomorphic statistical submersion. We have}
\[
\overline{\nabla}_{U}V=-F(\mathcal{H}\nabla_{U}^{\ast}(tV)),\hspace
{1.5cm}\overline{\nabla}_{U}^{\ast}V=-F(\mathcal{H}\nabla_{U}(tV)),
\]
\textit{moreover, if $\mathcal{H}\nabla_{U}P=0$, we get $\overline{\nabla}%
_{U}(FX)=F(\mathcal{H}\nabla_{U}^{\ast}X)$ and $\overline{\nabla}_{U}^{\ast
}(FX)=F(\mathcal{H}\nabla_{U}X)$. Furthermore, if $X$ is basic, then
$\overline{\nabla}_{U}(FX)=F(A_{X}^{\ast}U)$ and $\overline{\nabla}_{U}^{\ast
}(FX)=F(A_{X}U)$.}

\noindent\textbf{Lemma 4.7.}{\ \textit{Let $\pi:(\mathcal{M},\nabla
,g,J)\rightarrow(\mathcal{B},\widetilde{\nabla},\widetilde{g})$ be an
anti-invariant holomorphic statistical submersion. We have}}
\[
\mathcal{V}\nabla_{X}U=-F(\mathcal{H}\nabla_{X}^{\ast}(tU)),\hspace
{1.5cm}\mathcal{V}\nabla_{X}^{\ast}U=-F(\mathcal{H}\nabla_{X}(tU)),
\]
\textit{moreover, if $\mathcal{H}\nabla_{X}P=0$, we get $\mathcal{V}\nabla
_{X}(FY)=F(\mathcal{H}\nabla_{X}^{\ast}Y)$ and $\mathcal{V}\nabla_{X}^{\ast
}(FY)=F(\mathcal{H}\nabla_{X}Y)$.}

\vspace{0.2cm}The mean curvature vector of the affine connection $\nabla$ is
given by $N=\sum\varepsilon_{\alpha}T_{U_{\alpha}}U_{\alpha}$. If $\pi$ is an
anti-invariant holomorphic statistical submersion with conformal fiber, that
is, $T_{U}V=k\,g(U,V)$, then we find $k=\frac{\,N\,}{s}$, namely,
\begin{equation}
T_{U}V=\frac{\,1\,}{s}g(U,V)N \label{SN}%
\end{equation}
which yields from (\ref{BB}) that
\begin{equation}
T_{U}^{\ast}X=-\frac{1}{\,s\,}g(N,X)U. \label{SNN}%
\end{equation}
Changing $X$ to $tV$ in (\ref{SNN}), we get from Corollary 4.2
\begin{equation}
T_{U}V=\frac{1}{\,s\,}g(N,tV)tU. \label{SNNN}%
\end{equation}
Because of (\ref{SN}) and (\ref{SNNN}), we find $g(U,V)N=g(N,tV)tU$ which
yields from $PN=0$ that $sN=\sum\varepsilon_{\alpha}g(N,tU_{\alpha}%
)tU_{\alpha}=-tFN=N,$ that is, $(s-1)N=0.$ Thus we have

\noindent\textbf{Proposition 4.8.} \textit{Let }$\pi:(\mathcal{M}%
,\nabla,g,J)\rightarrow(\mathcal{B},\widetilde{\nabla},\widetilde{g}$%
$)$\textit{ be an anti-invariant holomorphic statistical submersion with
conformal fiber. If $\mathcal{H}\nabla_{U}P=0$, then we get}

\textit{i)}\ \ \textit{the dimension of each fiber is one, or}

\textit{ii)}\ \ $\pi$ \textit{is an anti-invariant holomorphic statistical
submersion with isometric fiber.}

%\section{Anti-invariant holomorphic statistical submersions satisfying the certain condition.}

\vspace{0.2cm}Let $\pi:(\mathcal{M},\nabla,g,J)\rightarrow(\mathcal{B}%
,\widetilde{\nabla},\widetilde{g}$$)$ be an anti-invariant holomorphic
statistical submersion which the curvature tensor with respect to the affine
connection $\nabla$ of the total space satisfies (\ref{*}) and $\mathcal{H}%
\nabla_{X}P=0$. Changing $Z$ to $PZ$ in (\ref{R3}), we find
$c\{g(tU,Y)g(PX,PZ)+2g(tU,X)g(PY,PZ)\}=0$ which means that $c=0$ or $P=0$.
Hence we have

\vspace{.2cm}

\noindent\textbf{Theorem 4.9.} \label{th0} \textit{Let }$\pi:(\mathcal{M}%
,\nabla,g,J)\rightarrow(\mathcal{B},\widetilde{\nabla},\widetilde{g}$%
$)$\textit{ be an anti-invariant holomorphic statistical submersion. If the
total space satisfies the condition $(\ref{*})$ and $\mathcal{H}\nabla_{X}%
P=0$, then}

\textit{i)}\ \ \textit{the total space is flat, or}

\textit{ii)}\ \ \textit{$P=0$.}

\vspace{.2cm}

Next, we discuss a holomorphic statistical submersion such that $P=0$. Then we
find $tF=-I,$ $fF=0,$ $tf=0$ and $f^{2}=-I-Ft$. From Lemma 3.2, we can get
\newline

\noindent\textbf{Lemma 4.10.} \textit{If }$\pi:(\mathcal{M},\nabla
,g,J)\rightarrow(\mathcal{B},\widetilde{\nabla},\widetilde{g}$$)$\textit{ is a
holomorphic statistical submersion satisfying $P=0$, then we have}
\begin{align}
&  \mathcal{H}\nabla_{U}(tV)+T_{U}(fV)=t(\overline{\nabla}_{U}^{\ast
}V),\label{P8}\\
&  T_{U}(tV)+\overline{\nabla}_{U}(fV)=F(T_{U}^{\ast}V)+f(\overline{\nabla
}_{U}^{\ast}V),\label{P9}\\
&  T_{U}(FX)=t(T_{U}^{\ast}X),\label{P10}\\
&  \overline{\nabla}_{U}(FX)=F(\mathcal{H}\nabla_{U}^{\ast}X)+f(T_{U}^{\ast
}X),\label{P11}\\
&  \mathcal{H}\nabla_{X}(tU)+A_{X}(fU)=t(\mathcal{V}\nabla_{X}^{\ast
}U),\label{P12}\\
&  A_{X}(tU)+\mathcal{V}\nabla_{X}(fU)=F(A_{X}^{\ast}U)+f(\mathcal{V}%
\nabla_{X}^{\ast}U),\label{P13}\\
&  A_{X}(FY)=t(A_{X}^{\ast}Y),\label{P14}\\
&  \mathcal{V}\nabla_{X}(FY)=F(\mathcal{H}\nabla_{X}^{\ast}Y)+f(A_{X}^{\ast
}Y). \label{P15}%
\end{align}
\textit{Furthermore, if $X$ is basic, then $\mathcal{H}\nabla_{U}^{\ast
}X=A_{X}^{\ast}U$.}

From (\ref{K36}), (\ref{K40}) and $(\nabla^{\ast})^{\ast}=\nabla,$ we find
\newline

\noindent\textbf{Lemma 4.11.} \textit{If }$\pi:(\mathcal{M},\nabla
,g,J)\rightarrow(\mathcal{B},\widetilde{\nabla},\widetilde{g}$$)$\textit{ is a
holomorphic statistical submersion satisfying $P=0$, then we have}
\[%
\begin{array}
[c]{ll}%
\mathcal{H}\nabla_{U}X=-t(\overline{\nabla}_{U}^{\ast}(FX)), & \hspace
{1.5cm}\mathcal{H}\nabla_{U}^{\ast}X=-t(\overline{\nabla}_{U}(FX))\\[0.2cm]%
\mathcal{H}\nabla_{X}Y=-t(\mathit{\mathcal{V}\nabla}_{U}^{\ast}(FY)), &
\hspace{1.5cm}\mathcal{H}\nabla_{X}^{\ast}Y=-t(\mathit{\mathcal{V}\nabla}%
_{U}(FY)).
\end{array}
\]

By virtue of (\ref{P8}), (\ref{P10}) and $(T^{\ast})^{\ast}=T$ we get \newline

\noindent\textbf{Lemma 4.12.} \textit{If }$\pi:(\mathcal{M},\nabla
,g,J)\rightarrow(\mathcal{B},\widetilde{\nabla},\widetilde{g}$$)$\textit{ is a
holomorphic statistical submersion satisfying $P=0$ and $\overline{\nabla}%
_{U}f=0$, then we have}
\[%
\begin{array}
[c]{ll}%
T_{U}V=-T_{U}(FtV)=-t(T_{U}^{\ast}(tV)), & \hspace{1cm}T_{U}^{\ast}%
V=-T_{U}^{\ast}(FtV)=-t(T_{U}(tV)),\\[0.2cm]%
T_{U}(fV)=T_{fV}U=0, & \hspace{1cm}T_{U}^{\ast}(fV)=T_{fV}^{\ast}U=0,\\[0.2cm]%
T_{U}X=-F(T_{U}^{\ast}(FX)), & \hspace{1cm}T_{U}^{\ast}X=-F(T_{U}%
(FX)),\\[0.2cm]%
f(T_{U}X)=0, & \hspace{1cm}f(T_{U}^{\ast}X)=0.
\end{array}
\]

From (\ref{P13}), (\ref{P14}) and $(A^{\ast})^{\ast}=A$, we obtain \newline

\noindent\textbf{Lemma 4.13.} \textit{If }$\pi:(\mathcal{M},\nabla
,g,J)\rightarrow(\mathcal{B},\widetilde{\nabla},\widetilde{g}$$)$\textit{ is a
holomorphic statistical submersion satisfying $P=0$ and $\mathcal{V}\nabla
_{X}f=0$, then we have}
\[%
\begin{array}
[c]{ll}%
A_{X}U=-t(A_{X}^{\ast}(tU))=t(A_{tU}X), & \hspace{1cm}A_{X}^{\ast}%
U=-t(A_{X}(tU))=t(A_{tU}^{\ast}X),\\[0.2cm]%
A_{X}(fU)=0, & \hspace{1cm}A_{X}^{\ast}(fU)=0,\\[0.2cm]%
A_{X}Y=-F(A_{X}^{\ast}(FY)), & \hspace{1cm}A_{X}^{\ast}Y=-F(A_{X}%
(FY)),\\[0.2cm]%
f(A_{X}Y)=0, & \hspace{1cm}f(A_{X}^{\ast}Y)=0.
\end{array}
\]

%\vspace{.2cm}

%{\bf Corollary 5.2.}\ \ {\it If $\pi:(M,\nabla,g,J)\to(B,\widehat{\nabla},g_B)$ is a holomorphic statistical submersion satisfying $P=0$, then we obtain}
%\begin{eqnarray}
%& & T^\ast_UV=-t\{T_U(tV)+\overline{\nabla}_U(fV)\},
%\\
%& & \overline{\nabla}^\ast_UV=-F\{{\cal H}\nabla_U(tV)+T_U(fV)\}-f\{T_U(tV)+\overline{\nabla}_U(fV)\},
%\\
%& & {\cal H}\nabla^\ast_UX=-t(\overline{\nabla}_U(FX)),
%\\
%& & T^\ast_UX=-F(T_U(FX))-f(\overline{\nabla}_U(FX)),
%\\
%& & A^\ast_XU=-t\{A_X(tU)+{\cal V}\nabla_X(fU)\},
%\\
%& & {\cal V}\nabla^\ast_XU=-F\{{\cal H}\nabla_X(tU)+A_X(fU)\}-f\{A_X(tU)+{\cal V}\nabla_X(fU)\},
%\\
%& & {\cal H}\nabla^\ast_XY=-t({\cal V}\nabla_X(FY)),
%\\
%& & A^\ast_XY=-F(A_X(FY))-f({\cal V}\nabla_X(FY)).
%\end{eqnarray}

\vspace{.2cm}

Let $\pi:(\mathcal{M},\nabla,g,J)\rightarrow(\mathcal{B},\widetilde{\nabla
},\widetilde{g}$$)$ be a holomorphic statistical submersion satisfying
(\ref{*}). We assume $\overline{\nabla}_{U}f=0$. Changing $W$ to $fW$ in
(\ref{R22}), we obtain from Lemma 4.12
\[
c\{g(tV,X)g(fU,fW)-g(tU,X)g(fV,fW)\}=0
\]
which means that $c=0$ or $g(tV,X)g(fU,fW)-g(tU,X)g(fV,fW)=0$. Thus we get
$||\,f\,||^{2}g(FX,U)=0.$ From $F\neq0$, we get $f=0$. Hence we have\vspace
{0.2cm}

\noindent\textbf{Theorem 4.14.} \label{th1} \textit{Let }$\pi:(\mathcal{M}%
,\nabla,g,J)\rightarrow(\mathcal{B},\widetilde{\nabla},\widetilde{g}$%
$)$\textit{ be a holomorphic statistical submersion satisfying $P=0$. If the
total space satisfies the condition $(\ref{*})$ and $\overline{\nabla}_{U}%
f=0$, then}

\textit{i)}\ \ \textit{ the total space is flat, or}

\textit{ii)}\ \ \textit{$f=0$.}

\bigskip

Because of Theorems 4.9 and 4.14, we have \newline

\noindent\textbf{Theorem 4.15.} \textit{Let }$\pi:(\mathcal{M},\nabla
,g,J)\rightarrow(\mathcal{B},\widetilde{\nabla},\widetilde{g}$$)$\textit{ be a
holomorphic statistical submersion which the total space satisfies the
condition $(\ref{*})$ with non-zero constant $c$. If $\mathcal{H}\nabla
_{X}P=0$ and $\overline{\nabla}_{U}f=0$, then $P=0$ is equivalent to $f=0$.}
\bigskip

Finally, we give an example of anti-invariant holomorphic statistical submersion.

\begin{example}
Let $\pi:(%
%TCIMACRO{\U{211d} }%
%BeginExpansion
\mathbb{R}
%EndExpansion
_{2}^{4},\nabla,g,J)\rightarrow(%
%TCIMACRO{\U{211d} }%
%BeginExpansion
\mathbb{R}
%EndExpansion
_{1}^{2},\widetilde{\nabla},\widetilde{g})$ be a holomorphic statistical
submersion given in Example\ref{EX2}. Then $\pi$ is an anti-invariant.
\end{example}

\end{document}